\documentclass[journal]{IEEEtran}
\newcommand{\beq}{\begin{equation}}
\newcommand{\eeq}{\end{equation}}
\newcommand{\bsq}{\begin{subequations}}
\newcommand{\esq}{\end{subequations}}
\newcommand{\bq}{\begin{eqnarray}}
\newcommand{\eq}{\end{eqnarray}}
\newcommand{\bqn}{\begin{eqnarray*}}
\newcommand{\eqn}{\end{eqnarray*}}
\usepackage{cite}
\usepackage{dsfont}
\usepackage{mathrsfs}
\usepackage{url}
\usepackage{empheq}
\usepackage{float}
\usepackage{dsfont}
\usepackage[hidelinks]{hyperref}

\usepackage{algorithm}
\usepackage{algorithmic}

\usepackage{amsmath,amssymb,amsfonts}
\allowdisplaybreaks[4]
\usepackage{graphicx}
\usepackage{textcomp}
\usepackage{amsmath}
\usepackage{enumerate}
\usepackage{epstopdf}
\usepackage{array}
\usepackage{booktabs}
\usepackage{subfigure}
\usepackage{multirow}
\usepackage{soul, color}
\usepackage[usenames,dvipsnames]{xcolor}
\usepackage[latin1]{inputenc}
\usepackage{cases}

\newtheorem{lemma}{Lemma}

\newtheorem{proposition}{Proposition}

\soulregister\cite7 
\soulregister\citep7 
\soulregister\citet7 
\soulregister\ref7 
\soulregister\pageref7

\usepackage{bbding}

\usepackage{nomencl}
\makenomenclature
\usepackage{etoolbox}
\renewcommand\nomgroup[1]{%
  \item[\bfseries
  \ifstrequal{#1}{A}{Acronyms}{%
  \ifstrequal{#1}{S}{Symbols}{%
  \ifstrequal{#1}{U}{Units}{}}}%
]}

\begin{document}

\title{Achieving Social Optimum and Budget Balance via a Joint Electricity-Carbon Pricing Mechanism}

\author{
Yue~Chen,~\IEEEmembership{Member,~IEEE},
Changhong~Zhao,~\IEEEmembership{Senior Member,~IEEE}
\thanks{Y. Chen is with the Department of Mechanical and Automation Engineering, the Chinese University of Hong Kong, HKSAR, China. (email: yuechen@mae.cuhk.edu.hk)}
\thanks{C. Zhao is with the Department of Information Engineering, the Chinese University of Hong Kong, HKSAR, China. (email: chzhao@ie.cuhk.edu.hk). The work of C. Zhao was supported by the Hong Kong Research Grants Council through General Research Fund No. 14212822.}}
\markboth{Journal of \LaTeX\ Class Files,~Vol.~XX, No.~X, Feb.~2019}%
{Shell \MakeLowercase{\textit{et al.}}: Bare Demo of IEEEtran.cls for IEEE Journals}

\maketitle

\begin{abstract}
Decarbonizing electric grids is a crucial global endeavor in the pursuit of carbon neutrality. Taking carbon emissions from generation into account when pricing electricity usage is an essential way to achieve this goal. However, such pricing is not trivial due to the requirements of an effective electricity market, such as maintaining budget balance, providing incentives to motivate participants to follow the dispatch schedule, and minimizing the impact on affected parties compared to when they were in the traditional electricity market. Although existing joint electricity-carbon pricing mechanisms have shown promising performance in reducing emissions in power networks, they can hardly meet all the requirements. This paper proposes a novel joint electricity-carbon pricing mechanism based on primal-dual optimality condition-enabled transformation. An algorithm for determining the critical market parameter is developed.  The proposed pricing mechanism is proven to possess all the desired properties, including budget balance, individual rationality, dispatch-following incentive compatibility, and truthful-bidding incentive compatibility. These properties ensure the proposed mechanism can incentivize market participants to achieve carbon-aware social optimum in a self-organized and sustainable way. Numerical experiments show the advantages of the proposed pricing mechanism compared to the existing marginal-based and carbon emission flow-based pricing mechanisms.
\end{abstract}


\IEEEpeerreviewmaketitle

\section*{Nomenclature}
\addcontentsline{toc}{section}{Nomenclature}

\subsection{Indices, Sets, and Symbols}
\begin{IEEEdescription}[\IEEEusemathlabelsep\IEEEsetlabelwidth{$Q(p,r^{\pm},w)$}]
\item[$\mathcal{I}, \mathcal{J}, \mathcal{L}$] Set of generators, loads, lines.
\item[$\mathcal{P}_i$] Set of variables constrained by \eqref{eq:market-proposed-3}.
\item[$\mathcal{D}_j$] Set of variables constrained by \eqref{eq:market-proposed-4}.
\item[$\mathcal{Q}$] Set of variables constrained by \eqref{eq:market-proposed-6}-\eqref{eq:market-proposed-8}.
\item[$\mathcal{S}_1, \mathcal{S}_2, \mathcal{S}_3$] Subsidy related to network congestion \eqref{eq:S1}, carbon tax \eqref{eq:S2} and market clearing \eqref{eq:S3}, respectively. $\Delta S=\mathcal{S}_2+\mathcal{S}_3$.
\end{IEEEdescription}

\subsection{Parameters}
\begin{IEEEdescription}[\IEEEusemathlabelsep\IEEEsetlabelwidth{$Q(p,r^{\pm},w)$}]
\item[$I, J$] Number of generators and loads.
\item[$\overline{P}_i$] Maximum power production of generator $i$.
\item[$c_i$] Cost coefficient of generator $i$.
\item[$e_i$] Carbon emission factor of generator $i$.
\item[$\overline{D}_j$] Maximum consumption of load $j$.
\item[$b_j$] Marginal utility of load $j$.
\item[$\kappa$] A coefficient turning carbon emissions into monetary values.
\item[$F_l$] Capacity limit of transmission line $l$.
\item[$\pi_{il},\pi_{jl}$] node-to-line transfer factors.
\item[$\delta$] A ratio for setting the carbon tax.
\item[$\bar{p},\bar{d}$] Optimal dispatch determined by the traditional electricity market.
\item[$p^*,d^*$] Optimal dispatch determined by the carbon-aware social welfare maximization problem.
\end{IEEEdescription}

\subsection{Decision Variables}
\begin{IEEEdescription}[\IEEEusemathlabelsep\IEEEsetlabelwidth{$Q(p,r^{\pm},w)$}]
\item[$p_{i}$] Power production of generator $i$.
\item[$d_j$] Demand of load $j$.
\item[$\lambda$, $\underline{\mu}_i, \overline{\mu}_i$] Dual variables of problem \eqref{eq:social-optimum}.
\item[$\underline{\phi}_j, \overline{\phi}_j$, $\overline{\chi}_l,\underline{\chi}_l$] Dual variables of problem \eqref{eq:social-optimum}.
\item[$\lambda_i^{\text{T1}}, \lambda_j^{\text{T1}}$] Joint electricity-carbon prices for generator $i$ and load $j$ under Method T1.
\item[$\lambda_i^{\text{T2}}, \lambda_j^{\text{T2}}$] Joint electricity-carbon prices for generator $i$ and load $j$ under Method T2.
\item[$\lambda_i^{\text{pro}}, \lambda_j^{\text{pro}}$] Joint electricity-carbon prices for generator $i$ and load $j$ under the proposed Method.
\item[$\sigma_i, \sigma_j$] Node carbon intensity for generator $i$ and load $j$.
\item[$\tau, \overline{\alpha}_l , \underline{\alpha}_l,\eta$] Dual variables of problem \eqref{eq:market-proposed}.
\item[$\text{LMP}_i,\text{LMP}_j$] Locational marginal prices under the CEF method.
\item[$\text{LMP}_i^{\text{tr}},\text{LMP}_j^{\text{tr}}$] Locational marginal prices in the traditional electricity market.
\end{IEEEdescription}

\section{Introduction}

\IEEEPARstart{T}{he} urgent global climate change calls for worldwide efforts to reduce carbon emissions. The electric sector is one of the main sources of carbon emissions, which requires particular attention \cite{chen2020low}. In the U.S. and European Union, while electricity accounts for only 20\% of the final energy consumption, over 40\% of all energy-related emissions are caused by electricity generation \cite{carbonreport}. Therefore, decarbonizing the electric sector is a critical step toward carbon neutrality.

There have been a vast literature studying the low-carbon operation of energy systems. Most of them were based on the carbon emission flow (CEF) method \cite{kang2015carbon}, which allocates the carbon responsibility in a power network along the power flows from generation to demand. References \cite{sun2024modeling}, \cite{cheng2018modeling}, and \cite{yang2023improved} extended the traditional CEF method to incorporate uncertainties, multi-energy systems, and prosumers, respectively. Although the CEF method is straightforward, its model is difficult to solve due to the nonlinearity caused by matrix inversion operations. To address such difficulty, a Bayesian inference regression-based data driven approach was developed to improve accuracy and computational efficiency \cite{wang2021optimal}. An adaptive regression method assisted with a carbon pattern dictionary was proposed to facilitate computation \cite{ma2024calculating}. A sparse neural network-aided method was developed to turn the CEF constraints into mixed integer linear constraints \cite{sang2023encoding}.  The CEF method has been widely applied in low-carbon power system planning \cite{wei2021carbon,wang2024two}, scheduling \cite{li2023low, wang2021optimal}, and optimal power flow \cite{chen2023carbon} problems.

In addition to carbon-aware operation under a centralized scheme, recent studies have looked into the pricing mechanism design to encourage low-carbon generation and consumption behaviors. Energy-carbon integrated pricing was proposed and used for electricity-gas system coordination \cite{cheng2019low}, electric vehicle management \cite{yuan2022low}, and peer-to-peer energy trading \cite{lu2022peer}. The above work derived the joint electricity-carbon prices based on the CEF model. Considering that the CEF method may not well capture the marginal contributions of agents in a power network \cite{xie2024real}, an alternative pricing mechanism is based on marginal analysis. Similar to the traditional locational marginal price (LMP), reference \cite{fang2022locational} proposed to add the carbon-related costs into the objective function of the optimal power flow (OPF) problem, and use the dual variable of the power balance condition as the carbon-aware LMP. This method was applied to integrated energy systems \cite{jiang2018optimal}. 
Reference \cite{ruiz2010analysis} further analyzed the difference of carbon intensities and prices at varying locations in a power network based on sensitivity analysis. An iterative algorithm was developed in \cite{park2023decarbonizing} to obtain the carbon-aware LMP while taking into account the response of elastic demand. A dynamic locational marginal emission quantification method was established in \cite{valenzuela2023dynamic}, which can better adapt to the variations of renewable and demand.

Despite the effectiveness of existing joint electricity-carbon pricing mechanisms in various applications, there is much less work exploring the properties of those pricing mechanisms theoretically. An electricity market should possess several important properties, including but not limited to the satisfaction of physical constraints, incentive compatibility, and individual rationality \cite{tsaousoglou2022market}. These properties are essential to ensure that the market can run in a self-organized and sustainable way \cite{chen2022pricing}. However, some counterexamples provided by recent studies \cite{cheng2019low,xie2024real} show that the existing joint electricity-carbon pricing mechanisms may violate the required properties.

In this paper, we fill this research gap by proposing a novel joint electricity-carbon pricing mechanism with proofs of all the desired properties. Our main contributions are two-fold:

1) \emph{Mechanism Design}. A joint electricity-carbon pricing mechanism is developed in this paper. First, we propose a market clearing problem based on the primal-dual optimality condition of the carbon-aware social welfare maximization problem. Next, by analyzing the structure of the Lagrangian dual of the market clearing problem, we derive the joint electricity-carbon prices for generators and loads, as well as the appropriate carbon tax rate. Furthermore, an algorithm to determine the critical market parameter is established with the help of the convex combination technique. The proposed market design is novel and has not been reported in existing literature.

2) \emph{Proof of Properties}. We prove that the proposed pricing mechanism possesses all the desired properties that are necessary for running an effective electricity market. These properties include budget balance, individual rationality, dispatch-following incentive compatibility, and truthful-bidding incentive compatibility. To the best of our knowledge, \emph{no} existing joint electricity-carbon pricing mechanisms is able to meet all of these properties simultaneously. Moreover, case studies show that the proposed pricing mechanism can achieve carbon-aware social optimum with minimum impact on the net profits of generators and loads.

The rest of this paper is organized as follows: The proposed joint electricity-carbon pricing mechanism is presented in Section \ref{sec-II}. In Section \ref{sec-III}, we prove several properties of the proposed pricing mechanism and compare it with the traditional pricing mechanisms. Numerical experiments are provided in Section \ref{sec-IV} with conclusions in Section \ref{sec-V}.

\section{Joint Electricity-Carbon Pricing Mechanism}
\label{sec-II}
We consider a power system with $I$ generators and $J$ loads, indexed by $i \in \mathcal{I}$ and $j \in \mathcal{J}$, respectively. To participate in the electricity market, each generator $i \in \mathcal{I}$ bids its maximum power production $\overline{P}_i$ and the marginal generation cost $c_i$; each load $j \in \mathcal{J}$ bids its maximum consumption quantity $\overline{D}_j$ and the marginal utility $b_j$. The carbon emission factor of each generator $i \in \mathcal{I}$ is denoted by $e_i$. 
Let $\kappa$ denote a coefficient that turns the carbon emission into a monetary value. In the following, we first introduce the social optimal problem that maximizes the carbon-aware social welfare, as well as two existing joint electricity-carbon pricing methods. Then, a new joint electricity-carbon pricing mechanism is proposed. 



\subsection{Carbon-Aware Social Welfare Maximization}
\label{sec-II-A}
To achieve a good balance between the economic costs and carbon emissions, the system operator solves the modified optimal power flow problem \eqref{eq:social-optimum} to determine the social optimal dispatch.
\bsq \label{eq:social-optimum}
\begin{align}
     & O^*:=\max_{p_i, \forall i, d_j,\forall j}~ \sum_{j\in \mathcal{J}} b_j d_j - \sum_{i \in \mathcal{I}} (c_i+\kappa e_i) p_i, \label{eq:social-optimum-1}\\
    & \mbox{s.t.}~ \sum_{i \in \mathcal{I}} p_i = \sum_{j \in \mathcal{J}} d_j: ~ \lambda, \label{eq:social-optimum-2}\\
    & 0 \le p_i \le \overline{P}_i,\forall i \in \mathcal{I}:~ \underline{\mu}_i, \overline{\mu}_i, \label{eq:social-optimum-3}\\
    & 0 \le d_j \le \overline{D}_j,\forall j \in \mathcal{J}:~\underline{\phi}_j, \overline{\phi}_j, \label{eq:social-optimum-4} \\
    & -F_l \le \sum_{i \in \mathcal{I}} \pi_{il} p_i - \sum_{j \in \mathcal{J}} \pi_{jl} d_j \le F_l,\forall l \in \mathcal{L}:~ \overline{\chi}_l,\underline{\chi}_l. \label{eq:social-optimum-5}
\end{align}
\esq
The objective function \eqref{eq:social-optimum-1} is to maximize the social welfare, i.e., the total consumer utility minus the total generation and carbon costs. Let $O^*$ denote the optimal objective value. Constraint \eqref{eq:social-optimum-2} ensures power balance. The dispatch of generators $p_i,\forall i$ and loads $d_j,\forall j$ are within their respective bounds   \eqref{eq:social-optimum-3} and \eqref{eq:social-optimum-4}. The power flow limits of transmission lines are given in \eqref{eq:social-optimum-5}, where $\pi_{il}$ and $\pi_{jl}$ are the node-to-line transfer factors. $\mathcal{L}$ is the set of transmission lines, each indexed by $l$. Dual variables associated with these constraints are denoted as $\lambda$, $\underline{\mu}_i, \overline{\mu}_i$, $\underline{\phi}_j, \overline{\phi}_j$, $\overline{\chi}_l,\underline{\chi}_l$. The value of a variable at the optimum of \eqref{eq:social-optimum} is indicated by the superscript $(\cdot)^*$. 

In existing literature, there are two mainstream methods for setting the joint electricity-carbon prices considering the power network constraints. Under both methods, the carbon tax rates are set as $\kappa$, i.e., each generator $i$ pays a carbon tax $\kappa e_i p_i$. In addition, each generator $i$ gets $\lambda_i^{(.)}p_i$ from the market and each load $j$ pays $\lambda_j^{(.)}d_j$ to the market; where $\lambda_i^{(.)}$ and $\lambda_j^{(.)}$ are the joint electricity-carbon prices depending on the pricing mechanisms. Specifically,

\begin{itemize}
    \item \emph{Marginal-Based Method} (\emph{T1 Method}) \cite{fang2022locational}. The operator solves the carbon-aware social welfare maximization problem \eqref{eq:social-optimum} and derives the locational marginal prices following a similar way as the traditional electricity market does. As a result, the prices are
    \bsq \label{eq:price-T1}
    \begin{align}
        \text{Generator}~i: ~ & \lambda_i^{\text{T1}}=\lambda^*+\sum_{l \in \mathcal{L}} (\underline{\chi}_l^*-\overline{\chi}_l^*) \pi_{il}, \\
        \text{Load}~j:~& \lambda_j^{\text{T1}}=\lambda^*+\sum_{l \in \mathcal{L}} (\underline{\chi}_l^*-\overline{\chi}_l^*) \pi_{jl}.
    \end{align}
    \esq
    \item \emph{Carbon Emission Flow-Based Method} (\emph{T2 Method}) \cite{cheng2019low}. Under this method, a bilevel optimization model is formulated. In the upper level, the operator solves \eqref{eq:opf} to derive the locational marginal prices $\text{LMP}_i$, $i\in\mathcal{I}$ and $\text{LMP}_j$, $j \in \mathcal{J}$.
    \bsq \label{eq:opf}
    \begin{align} \label{eq:opf-obj}
        \min_{p_i,\forall i\in\mathcal{I}}~ & \sum_{i \in \mathcal{I}} c_i p_i, \\
        \mbox{s.t.}~ & \eqref{eq:social-optimum-2}, \eqref{eq:social-optimum-3}, \eqref{eq:social-optimum-5},\\
        ~ & \{d_j,\forall j\in\mathcal{J}\}~\text{is the optimal of \eqref{eq:load}}.
    \end{align}
    \esq
    Then, the operator calculates the node carbon intensity (NCI) for generators $\sigma_i,\forall i\in\mathcal{I}$ and loads $\sigma_j,\forall j\in\mathcal{J}$ using the carbon emission flow method \cite{kang2015carbon}. The weighted sum of the electricity price and NCI is set as the joint electricity-carbon price, i.e.,
    \bsq
    \begin{align}
        \text{Generator}~i: ~ & \lambda_i^{\text{T2}}=\text{LMP}_i + \kappa \sigma_i,\\
        \text{Load}~j:~& \lambda_j^{\text{T2}}=\text{LMP}_j + \kappa \sigma_j.
    \end{align}
    \esq
In the lower-level, the optimal consumption quantities of loads are determined by solving
\begin{align}\label{eq:load}
    \max_{d_j,\forall j\in\mathcal{J}}~ \sum_{j\in\mathcal{J}} (b_j  -  \lambda_j^{\text{T2}})d_j~~~~ \mbox{s.t.}~ \eqref{eq:social-optimum-4}
\end{align}
An iterative algorithm was developed in \cite{cheng2019low} to solve the bilevel model to obtain the $\lambda_i^{\text{T2}}$ and $\lambda_j^{\text{T2}}$ at equilibrium.
\end{itemize}

In the following, we first introduce the proposed joint electricity-carbon pricing mechanism. Then, the limitations of the existing pricing mechanisms and the desired properties of the proposed mechanism are discussed in Section \ref{sec-III}.

\subsection{The Proposed Pricing Mechanism}
In the paper, we propose a new primal-dual optimality condition-based pricing mechanism to determine the joint electricity-carbon prices. Specifically, first, we write down the dual problem of the carbon-aware social welfare maximization problem \eqref{eq:social-optimum}, which is
\begin{align}
    & \min_{\lambda, \underline{\mu}, \overline{\mu},\underline{\phi}, \overline{\phi},\underline{\chi},\overline{\chi}}~ \sum_{i \in \mathcal{I}} \overline{\mu}_i\overline{P}_i + \sum_{j \in \mathcal{J}} \overline{\phi}_j\overline{D}_j + \sum_{l \in \mathcal{L}} \overline{\chi}_l F_l + \sum_{l \in \mathcal{L}} \underline{\chi}_l F_l, \nonumber\\
    & \mbox{s.t.}~ -\lambda -\underline{\mu}_i + \overline{\mu}_i - \sum_{l \in \mathcal{L}} \pi_{il}\overline{\chi}_l + \sum_{l \in \mathcal{L}} \pi_{il}\underline{\chi}_l = -(c_i+\kappa e_i),\forall i, \nonumber\\
    & \lambda -\underline{\phi}_j +\overline{\phi}_j + \sum_{l \in \mathcal{L}} \pi_{jl} \overline{\chi}_l - \sum_{l \in \mathcal{L}} \pi_{jl} \underline{\chi}_l = b_j,\forall j, \nonumber\\
    & \underline{\mu}, \overline{\mu}, \underline{\phi}, \overline{\phi}, \underline{\chi}, \overline{\chi} \ge 0.
\end{align}

Instead of using \eqref{eq:social-optimum} or the OPF problem to clear market, we formulate the following market-clearing problem:
\bsq \label{eq:market-proposed}
\begin{align}
    & \max ~ \sum_{j \in \mathcal{J}} b_j d_j -\sum_{i \in \mathcal{I}} (c_i + \delta \kappa e_i) p_i, \label{eq:market-proposed-1}\\
    & \mbox{s.t.}~ \sum_{i \in \mathcal{I}} p_i =\sum_{j \in \mathcal{J}} d_j: ~ \tau, \label{eq:market-proposed-2}\\
     & 0 \le p_i \le \overline{P}_i,\forall i, \label{eq:market-proposed-3}\\
     & 0 \le d_j \le \overline{D}_j,\forall j, \label{eq:market-proposed-4}\\
     & -F_l \le \sum_{i \in \mathcal{I}} \pi_{il} p_i - \sum_{j \in \mathcal{J}} \pi_{jl} d_j \le F_l,\forall l :~\overline{\alpha}_l,\underline{\alpha}_l, \label{eq:market-proposed-5}\\
    & -\lambda - \underline{\mu}_i + \overline{\mu}_i -\sum_{l \in \mathcal{L}} \pi_{il}(\overline{\chi}_l-\underline{\chi}_l) = -(c_i+\kappa e_i),\forall i, \label{eq:market-proposed-6}\\
    & \lambda -\underline{\phi}_j +\overline{\phi}_j + \sum_{l \in \mathcal{L}} \pi_{jl} (\overline{\chi}_l - \underline{\chi}_l) = b_j,\forall j, \label{eq:market-proposed-7}\\
    & \underline{\mu}, \overline{\mu}, \underline{\phi}, \overline{\phi}, \underline{\chi}, \overline{\chi} \ge 0, \label{eq:market-proposed-8}\\
    &  \sum_{i \in \mathcal{I}} \overline{\mu}_i \overline{P}_i + \sum_{j \in \mathcal{J}} \overline{\phi}_j \overline{D}_j + \sum_{l \in \mathcal{L}} \overline{\chi}_l F_l + \sum_{l \in \mathcal{L}} \underline{\chi}_l F_l \nonumber\\
    & \le \sum_{j \in \mathcal{J}} b_j d_j -\sum_{i \in \mathcal{I}} (c_i+ \kappa e_i)p_i:~\eta. \label{eq:market-proposed-9}
\end{align}
\esq
In the objective function \eqref{eq:market-proposed-1}, we use a factor $\delta \in [0,1]$ to scale the carbon emission costs. We will see later that  $\delta \kappa$ can be interpreted as the carbon tax rate under the proposed pricing mechanism. Constraints
\eqref{eq:market-proposed-2}-\eqref{eq:market-proposed-5} are the primal constraints of \eqref{eq:social-optimum} while \eqref{eq:market-proposed-6}-\eqref{eq:market-proposed-8} are the dual constraints. Inequality \eqref{eq:market-proposed-9} together with the weak duality condition yields equality for \eqref{eq:market-proposed-9} \cite{boyd2004convex}. Hence, \eqref{eq:market-proposed-2}-\eqref{eq:market-proposed-9} is the primal-dual optimality condition of the social optimal problem \eqref{eq:social-optimum}. As a result, problems \eqref{eq:social-optimum} and \eqref{eq:market-proposed} have the same optimal dispatch $(p^*,d^*)$. $\tau$, $\overline{\alpha}_l,\underline{\alpha}_l$, and $\eta$ denote the dual variables of \eqref{eq:market-proposed-2}, \eqref{eq:market-proposed-5}, \eqref{eq:market-proposed-9}, respectively.

Define the set of variables constrained by \eqref{eq:market-proposed-3} as $\mathcal{P}_i,\forall i\in \mathcal{I}$, by \eqref{eq:market-proposed-4} as $\mathcal{D}_j,\forall j\in\mathcal{J}$, and by \eqref{eq:market-proposed-6}-\eqref{eq:market-proposed-8} as $\mathcal{Q}$. The Lagrangian dual of \eqref{eq:market-proposed} is given by:
\begin{align}\label{eq:Lagrangian-2}
    \min_{\tau, \underline{\alpha} \ge 0, \overline{\alpha} \ge 0, \eta \ge 0}~ \max_{p_i \in \mathcal{P}_i, d_j \in \mathcal{D}_j, (\lambda, \overline{\mu},\underline{\mu}, \overline{\phi},\underline{\phi},\underline{\chi},\overline{\chi}) \in \mathcal{Q}} \mathscr{L},
\end{align}
where
\begin{align}
    \mathscr{L}=~ &  \sum_{j \in \mathcal{J}} b_j d_j -\sum_{i \in \mathcal{I}} (c_i + \delta \kappa e_i) p_i + \tau (\sum_{i \in \mathcal{I}}  p_i-\sum_{j \in \mathcal{J}} d_j ) \nonumber\\
    ~ & + \sum_{l \in \mathcal{L}} \underline{\alpha}_l \left(F_l + \sum_{i \in \mathcal{I}} \pi_{il}p_i - \sum_{j \in \mathcal{J}} \pi_{jl}d_j \right) \nonumber\\
    ~ &  +\sum_{l \in \mathcal{L}} \overline{\alpha}_l \left(F_l - \sum_{i \in \mathcal{I}} \pi_{il}p_i +  \sum_{j \in \mathcal{J}} \pi_{jl}d_j\right) \nonumber\\
    ~ &  +\eta \left(-\sum_{i \in \mathcal{I}} \overline{\mu}_i \overline{P}_i - \sum_{j \in \mathcal{J}} \overline{\phi}_j \overline{D}_j - \sum_{l \in \mathcal{L}} \overline{\chi}_l F_l \right. \nonumber\\
    ~ & \left.-\sum_{l \in \mathcal{L}} \underline{\chi}_l F_l + \sum_{j \in \mathcal{J}} b_j d_j -\sum_{i \in \mathcal{I}} (c_i+ \kappa e_i)p_i\right).
\end{align}

Solving problem \eqref{eq:market-proposed}, we can obtain the $\tau^*$, $\eta^*$, $\underline{\alpha}_l^*, \overline{\alpha}_l^*,\forall l$ at the optimum. Fixing these values, the inner $\max\mathscr{L}$ can be separated into \eqref{eq:strategic-gen} and \eqref{eq:strategic-load} for each agent (generator or load), whose optimal solutions are $p^*$ and $d^*$, respectively.

Particularly, each generator $i \in \mathcal{I}$ solves
\begin{align} \label{eq:strategic-gen}
    \max_{p_i \in \mathcal{P}_i}~ -(c_i+\delta \kappa e_i) p_i + \lambda_i^{\text{pro}}p_i, 
\end{align}
and 
each load $j\in \mathcal{J}$ solves
\begin{align} \label{eq:strategic-load}
\max_{d_j \in \mathcal{D}_j} b_j d_j -\lambda_j^{\text{pro}}d_j,
\end{align}
where
\bsq
\begin{align}
    \lambda_i^{\text{pro}}=& ~\tau^* + \sum_{l \in \mathcal{L}} (\underline{\alpha}_l^*- \overline{\alpha}_l^* )\pi_{il}  - \eta^* (c_i+\kappa e_i), \\
\lambda_j^{\text{pro}}=& ~\tau^*+\sum_{l\in \mathcal{L}} (\underline{\alpha}_l^* - \overline{\alpha}_l^*) \pi_{jl}-\eta^* b_j. 
\end{align}
\esq

\begin{figure}[ht]
  \centering
\includegraphics[width=0.4\textwidth]{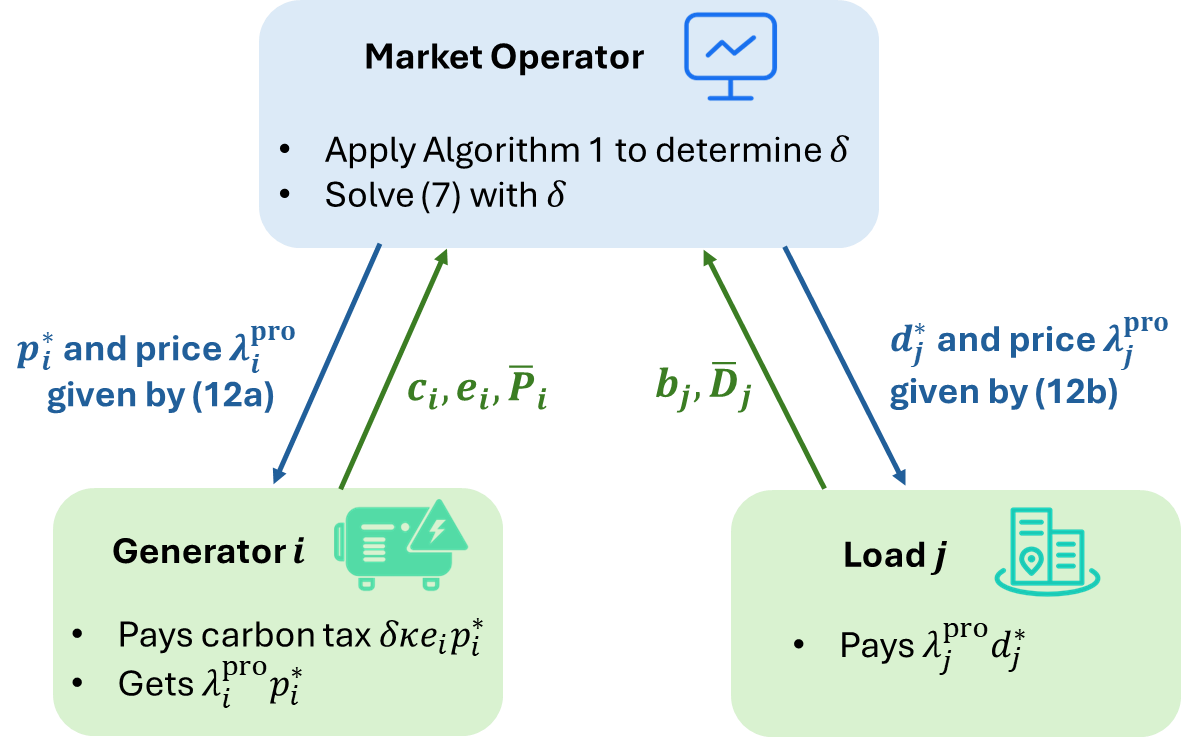}\\
  \caption{Overall procedure of the proposed pricing mechanism.}\label{fig:procedure}
\end{figure}

Based on the above facts, the proposed joint electricity-carbon pricing mechanism is designed as follows: Each generator $i$ pays a carbon tax with the rate of $\delta \kappa$ and gets paid at $\lambda_i^{\text{pro}}$; and each load $j$ pays at $\lambda_j^{\text{pro}}$. It is worth noting that the value of $\eta$ at the optimum of \eqref{eq:market-proposed} may not be unique. If multiple solutions occur, we use the one with the smallest value as $\eta^*$ to reduce the impact of carbon-related prices on the affected agents. The overall procedure of the proposed pricing mechanism is shown in Fig. \ref{fig:procedure}.

Under the proposed pricing mechanism, the subsidy required to support the market operation can be calculated by
\begin{align}\label{eq:subsidy}
    \mathcal{S} = ~ & \text{generator revenue}-\text{carbon tax}-\text{load payment}, \nonumber\\
    =~ & \sum_{i\in\mathcal{I}} \lambda_i^{\text{pro}} p_i^* - \sum_{i\in \mathcal{I}} \delta \kappa e_i p_i^* - \sum_{j\in \mathcal{J}} \lambda_j^{\text{pro}} d_j^*, \nonumber\\
    =~ & \mathcal{S}_1 + \mathcal{S}_2 +\mathcal{S}_3,
\end{align}
where
\bsq
\begin{align}
    \mathcal{S}_1 =~ & \sum_{l\in\mathcal{L}} (\underline{\alpha}_l^*-\overline{\alpha}_l^*) \left(\sum_{i\in\mathcal{I}} \pi_{il}p_i^* - \sum_{j\in \mathcal{J}} \pi_{jl} d_j^*\right), \label{eq:S1}\\
    \mathcal{S}_2=~ & -\delta \kappa \sum_{i\in\mathcal{I}} e_i p_i^*, \label{eq:S2}\\
    \mathcal{S}_3 =~ & \eta^* \left(-\sum_{i\in\mathcal{I}} (c_i+\kappa e_i) p_i^*
 + \sum_{j\in\mathcal{J}} b_j d_j^*\right). \label{eq:S3}
\end{align}
\esq
The first term $\mathcal{S}_1 = -\sum_{l\in\mathcal{L}} (\underline{\alpha}_l^*+\overline{\alpha}_l^*) F_l\le 0$ since $\underline{\alpha}_l^*, \overline{\alpha}_l^*, F_l$ are all nonnegative; and it is related to network congestion, meaning that the market can gain extra revenue when the network is congested. This revenue is called ``financial transmission right'' and is usually used for power network reinforcement \cite{sarkar2008comprehensive}. The second term $\mathcal{S}_2 \le 0$ is the carbon tax revenue of the government. The third term $\mathcal{S}_3 \ge 0$ since $\eta^* \ge 0$ and the optimal value of \eqref{eq:social-optimum} is $O^*\geq 0$. Define
\begin{eqnarray}
\Delta S:=\mathcal{S}_2+\mathcal{S}_3.\nonumber
\end{eqnarray} 
In the following, we prove that the proposed joint electricity-carbon pricing mechanism achieves budget balance.


\section{Properties of the Pricing Mechanism}
\label{sec-III}
In this section, we first discuss the desired properties of a pricing mechanism for a joint electricity-carbon market and the limitations of existing mechanisms. Then, we prove that the proposed mechanism possesses all the desired properties.

\subsection{Desired Properties of a Pricing Mechanism}
\label{sec-III-A}
The joint electricity-carbon pricing mechanism should preserve the desired properties of the traditional electricity markets. Specifically, in a traditional electricity market, the social welfare maximization problem \eqref{eq:social-optimum} with $\kappa=0$ is solved to determine the locational marginal prices $\text{LMP}_i^{\text{tr}},\forall i\in\mathcal{I}$, $\text{LMP}_j^{\text{tr}},\forall j\in\mathcal{J}$ and the optimal dispatch $(\bar p, \bar d)$. The main properties of a traditional electricity market and their counterparts for a joint electricity-carbon pricing mechanism are described below and summarized in Table \ref{tab:property}.

\begin{table*}[t]
        \renewcommand{\arraystretch}{1.8}
        \renewcommand{\tabcolsep}{0.5em}
        \centering
        \caption{Properties of Traditional Electricity Markets and Counterparts for a Joint Electricity-Carbon Pricing Mechanism}
        \label{tab:property}
        \begin{tabular}{m{0.15\linewidth}<{\centering}m{0.3\linewidth}<{\centering}m{0.32\linewidth}<{\centering}}
         \toprule
         & Traditional Electricity Market &   Joint Electricity-Carbon Pricing Mechanism\\
         \midrule 
      Budget Balance & When network is not congested, total generator revenue equals total load payment. $\sum_{j\in\mathcal{J}} \text{LMP}_j^{\text{tr}}\bar{d}_j$=$\sum_{i\in\mathcal{I}} \text{LMP}_i^{\text{tr}}\bar{p}_i$ &  When network is not congested, total generator revenue equals total load and carbon payment.
      \\
     Individual Rationality  & The net profit of each generator and the net utility of each load are nonnegative. $\text{LMP}_i^{\text{tr}} \bar p_i - c_i\bar p_i \ge 0,\forall i\in\mathcal{I}$, $b_j \bar d_j - \text{LMP}_j^{\text{tr}}\bar{d}_j \ge 0,\forall j\in \mathcal{J}$ &  The net profit of each generator and the net utility of each load are nonnegative. $\text{price}_i \times p_i-(c_i+ \text{carbon tax rate} \times e_i) p_i\ge 0,\forall i$, 
 $b_j d_j -\text{price}_j \times d_j \ge 0,\forall j$\\
    Dispatch-Following Incentive Compatibility & Given the prices $\text{LMP}_i^{\text{tr}},\forall i$ and $\text{LMP}_j^{\text{tr}},\forall j$, the problems \eqref{eq:trLMP-1} and \eqref{eq:trLMP-2} are solved at the optimal solution of \eqref{eq:social-optimum} with $\kappa=0$. & Given the joint electricity-carbon prices, the problems \eqref{eq:req-1} and \eqref{eq:req-2} are solved at the optimal solution of \eqref{eq:social-optimum} ($\kappa \ne 0$).\\
    Truthful-Bidding Incentive Compatibility & It is optimal for every agent to reveal its true cost function and parameters. & It is optimal for every agent to reveal its true cost function and parameters.\\
     \bottomrule
        \end{tabular}
\end{table*}

\begin{enumerate}
    \item \emph{Budget Balance}. In traditional electricity markets, when the power network is not congested, the total load payment $\sum_{j\in\mathcal{J}} \text{LMP}_j^{\text{tr}}\bar{d}_j$ equals the total generator revenue $\sum_{i\in\mathcal{I}} \text{LMP}_i^{\text{tr}}\bar{p}_i$. In other words, budget balance is attained, implying that the operator does not require additional subsidy to run the market, nor does it make a profit from the market. Correspondingly, under the joint electricity-carbon pricing mechanism, when the power network is not congested, we should have that the total carbon tax payment plus the total load payment equals the total generation revenue.
    \item \emph{Individual Rationality}. The net profit of each generator and the net utility of each load are both nonnegative, i.e., $\text{LMP}_i^{\text{tr}} \bar p_i - c_i\bar p_i \ge 0,\forall i\in\mathcal{I}$ and $b_j \bar d_j - \text{LMP}_j^{\text{tr}}\bar{d}_j \ge 0,\forall j\in \mathcal{J}$.
    Intuitively, all agents are willing to participate in the market since they get better off. Under the joint electricity-carbon pricing mechanism, we should have the same property considering the carbon tax payments. 
    \item \emph{Dispatch-Following Incentive Compatibility}. In traditional electricity markets, given the LMPs,  the individually optimal dispatches of generators and loads together achieve social optimum. Specifically, for each $i\in \mathcal{I}$ and $j\in \mathcal{J}$, the following optimization problems: 
    \bsq
    \begin{align}
        \text{Generator}~i:&~ \max_{p_i \in [0,\overline{P}_i]}~ - c_i p_i +\text{LMP}_i^{\text{tr}} p_i \label{eq:trLMP-1}\\
        \text{Load}~j:&~ \max_{d_j \in [0,\overline{D}_j]}~ b_j d_j -\text{LMP}_j^{\text{tr}} d_j \label{eq:trLMP-2}
    \end{align}
    \esq
    are solved at the social optimal solution $(\bar p,\bar d)$ for \eqref{eq:social-optimum} with $\kappa=0$, i.e., without considering carbon tax. In the joint electricity-carbon pricing mechanism, such relation should also hold between individual dispatch \eqref{eq:req-1}, \eqref{eq:req-2} and social welfare maximization \eqref{eq:social-optimum} with $\kappa\neq 0$.
    \bsq
    \begin{align}
        \text{Generator}~i:~ \max_{p_i \in [0,\overline{P}_i]}~ & - (c_i +\text{carbon tax rate}) \times p_i \nonumber\\
        ~ & +\text{price}_i \times p_i \label{eq:req-1}\\
        \text{Load}~j:~ \max_{d_j \in [0,\overline{D}_j]}~ & b_j d_j -\text{price}_j d_j \label{eq:req-2}
    \end{align}
    \esq
\item \emph{Truthful-Bidding Incentive Compatibility}. It is optimal for every agent to reveal its true cost function and operational parameters such as the marginal generation cost $c_i$, the carbon emission factor $e_i$, the marginal utility $b_j$, and the maximum quantities $\overline{P}_i, \overline{D}_j$. 
\end{enumerate}

\subsection{Limitations of Existing Pricing Mechanisms}
Two existing pricing mechanisms are introduced in Section \ref{sec-II-A}, namely the marginal-based method (T1 method) and the CEF-based method (T2 method). In the following, we show that the existing mechanisms can hardly satisfy all the required properties mentioned above.


\subsubsection{Marginal-Based Method (T1 Method)} 
First, we show that this method cannot achieve budget balance. Under the T1 method, the prices are given by \eqref{eq:price-T1}. Specially, when the power network is not congested, we have $\lambda_i^{\text{T1}}=\lambda^*=\lambda_j^{\text{T1}}$, $\forall i\in \mathcal{I}$, $j\in \mathcal{J}$. Therefore, the total generation revenue minus the total load and carbon payment is
\begin{align}
    \sum_{i\in\mathcal{I}} \kappa e_i p_i^* >0,
\end{align}
i.e., budget balance is violated and the operator makes an extra profit from the market. This extra profit is due to the significant increase in the price of loads, as proven in the Appendix A of \cite{cheng2019low}. In other words, the T1 method achieves carbon-aware social optimum at the cost of consumers, making them reluctant to accept this pricing mechanism.


\subsubsection{CEF-Based Method (T2 Method)} According to the simulations in Section \ref{sec-IV}, we can easily find that budget balance and dispatch-following incentive compatibility properties do not hold under the CEF-based method. Particularly, as shown in the 7th column of Table \ref{tab:Proposed-T1-T2}, the subsidy under T2 is negative. This means that, similar to the T1 method, the operator makes an extra profit from the market. Moreover, as shown in the last column, the social welfare under T2 is not the optimal value of \eqref{eq:social-optimum}, indicating that the prices under CEF-based method fail to encourage strategic agents to follow the carbon-aware social optimal dispatch strategies.


\subsection{Properties of the Proposed Pricing Mechanism}
While the existing pricing mechanisms may violate some of the desired properties in Section \ref{sec-III-A}, in the following, we prove that the proposed joint electricity-carbon pricing mechanism has all those properties. 

\begin{proposition}(\emph{Budget Balance}) \label{prop-1}
There exists $\delta \in [0,1)$ such that $\Delta \mathcal{S}=0$. At such $\delta$, if there is no power network congestion, then $\mathcal{S}=0$.
\end{proposition}

\begin{proposition}(\emph{Individual Rationality})\label{prop-2}
    Under the proposed pricing mechanism, the net profit of each generator and the net utility of each load are both nonnegative, i.e., denoting the optimal objective values of \eqref{eq:strategic-gen} and \eqref{eq:strategic-load} by $\text{NP}_{i}^*$ and $\text{NU}_j^*$, respectively, we have $\text{NP}_{i}^*\ge 0,\forall i\in \mathcal{I}$ and $\text{NU}_j^* \ge 0,\forall j \in \mathcal{J}$.
\end{proposition}

\begin{proposition}(\emph{Dispatch-Following Incentive Compatibility})\label{prop-3}
    Under the proposed pricing mechanism, the strategic generators and loads would produce the socially optimal dispatch $(p^*,d^*)$ by themselves in ex-post, given the realized joint electricity-carbon prices $\lambda_i^{\text{pro}}, \lambda_j^{\text{pro}}$ and carbon tax rate $\delta \kappa$.
\end{proposition}

\begin{proposition}(\emph{Truthful-Bidding Incentive Compatibility})\label{prop-4}
    Under the proposed pricing mechanism, it is optimal for price-taking agents to bid their true cost functions and true operational parameters.
\end{proposition}

The proofs of Propositions \ref{prop-1}-\ref{prop-4} above can be found in Appendices \ref{appendix-1}-\ref{appendix-4}, respectively. These propositions show that the proposed pricing mechanism possesses all the desired properties and is hence promising. Furthermore, according to Proposition \ref{prop-1}, the parameter $\delta$ should be properly set to ensure budget balance. Here, we propose Algorithm \ref{algo} to determine the value of $\delta$. It is based on the fact that a convex combination of the optimal solutions $\upsilon_{0}^*$ and $\upsilon_{\tilde{\delta}}^*$ under $\delta=0$ and $\delta = \tilde{\delta}$, respectively, is an optimal solution for the corresponding $\delta$. For more details, please refer to Lemma \ref{lemma-2} in Appendix \ref{appendix-1}.

\begin{algorithm}[htbp]
\caption{Determination of Parameter $\delta$}\label{algo}
\begin{algorithmic}[1]
\STATE{Starting from $\delta=1$, keep decreasing $\delta$ until the $\eta_{\delta}^*$ obtained by \eqref{eq:market-proposed} becomes non-zero. Denote the last (lowest) $\delta$ with $\eta_{\delta}^*=0$ by $\tilde{\delta}$.}
\STATE{Solve the problem \eqref{eq:market-proposed} with $\delta=0$ and obtain the corresponding optimal solutions $p_i^*|_{\delta=0}$, $d_j^*|_{\delta=0}$, $\eta^*_{0}$.}
\STATE{Solve the problem \eqref{eq:market-proposed} with $\delta=\tilde{\delta}$ and obtain the corresponding optimal solutions $p_i^*|_{\delta=\tilde{\delta}}$, $d_j^*|_{\delta=\tilde{\delta}}$, $\eta^*_{\tilde{\delta}}$.}
\STATE{Solve the following equation:
\begin{footnotesize}
\begin{align} \label{eq:qua}
  0 = &  -x \tilde{\delta} \kappa \sum_i e_i \left[ x p_i^*|_{\delta=\tilde{\delta}} + (1-x) p_i^*|_{\delta=0}\right] \nonumber\\
  ~ &  + (1-x) \eta^*_0 \left[ -\sum_i (c_i+\kappa e_i) \left(x p_i^*|_{\delta=\tilde{\delta}} + (1-x) p_i^*|_{\delta=0}\right) \right. \nonumber\\
  ~ &  \left.+ \sum_j b_j \left(x d_j^*|_{\delta=\tilde{\delta}} + (1-x) d_j^*|_{\delta=0}\right)\right].
\end{align}
\end{footnotesize}
Denote the solution by $x^*$.}
\STATE{Let $\delta = x^* \tilde{\delta}$.}
\end{algorithmic}
\end{algorithm}


\section{Numerical Experiments}
\label{sec-IV}
We first use a simple system to validate the proposed pricing mechanism and compare it to the traditional ones. Then, larger systems are tested to further analyze the performance and scalability of the proposed pricing mechanism. The system data can be found in \cite{Data}. The time slot considered is 1 hour.

\subsection{Simple System}
\label{sec:simple system}

We first test the performance of the proposed pricing mechanism using a simple system with 6 generators and 8 loads. The line capacity $F_l,\forall l\in \mathcal{L}$ are set large enough so that the power network is not congested, i.e., \eqref{eq:social-optimum-5} is not binding. The other parameters are given in Tables \ref{tab:para-gen}-\ref{tab:para-load}. $\kappa=$0.07\$/kgCO$_2$. 

\begin{table}[t]
\renewcommand{\arraystretch}{1.3}
\renewcommand{\tabcolsep}{0.5em}
\centering
\caption{Parameters of generators}
\label{tab:para-gen}
\begin{tabular}{ccccccc} 
\hline
Generator & 1 & 2 & 3 & 4 & 5 & 6 \\
\hline
$c_i$ (\$/kWh) & 0.472 & 0.480 & 0.502 & 0.473 & 0.492 & 0.512 \\
$e_i$ (kgCO$_2$/kWh) & 0.9 & 0.8 & 0.8 & 0.2 & 0.3 & 0.3 \\
$\overline{P}_i$ (kW) &  800 & 800 & 500 & 550 & 300 & 400 \\
\hline
\end{tabular}
\end{table}

\begin{table}[t]
\renewcommand{\arraystretch}{1.3}
\renewcommand{\tabcolsep}{0.5em}
\centering
\caption{Parameters of loads}
\label{tab:para-load}
\begin{tabular}{ccccccccc} 
\hline
Load & 1 & 2 & 3 & 4 & 5 & 6 & 7 & 8 \\
\hline
$b_j$ (\$/kWh) & 0.78 & 0.78 & 0.85 & 0.67 & 0.85 & 0.73 & 0.75 & 0.84 \\
$\overline{D}_i$ (kW) &  350 & 340 & 420 & 500 & 200 & 330 & 280 & 250 \\
\hline
\end{tabular}
\end{table}

\subsubsection{Effectiveness}
\begin{figure}[ht]
  \centering
\includegraphics[width=0.48\textwidth]{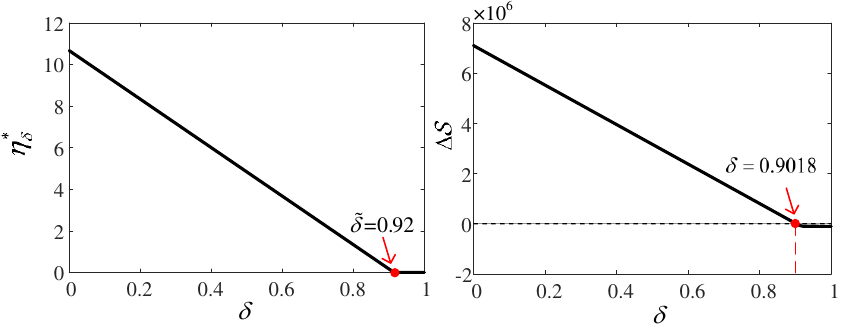}\\
  \caption{$\eta_{\delta}^*$ (\$/kWh) and $\Delta \mathcal{S}$ (\$) under different $\delta$.}\label{fig:delta}
\end{figure}

We apply Algorithm \ref{algo} to determine the value of parameter $\delta$. Fig. \ref{fig:delta} shows how $\eta_{\delta}^*$ and $\Delta S$ change with $\delta$. As shown, $\tilde \delta=0.92$ and when $\delta \in (0,\tilde{\delta})$, $\eta_{\delta}^*$ decreases linearly with $\delta$, which validates Lemma \ref{lemma-2} in Appendix \ref{appendix-1}. We can determine $\delta=0.9018$ that makes $\Delta S=0$ both by Algorithm \ref{algo} and the figure. Under this $\delta$, we test the proposed pricing mechanism. 
The total generator revenue is \$1,421,658, the total carbon tax payment is \$96,962, and the total load payment is \$1,324,696. Budget balance is achieved since the total revenue equals the total payment, which verifies Proposition \ref{prop-1}. Moreover, the carbon-aware social welfare under the proposed market is \$665,830, which is the same as the optimal objective value of \eqref{eq:social-optimum}. This implies that the proposed pricing mechanism can achieve carbon-aware social optimum. The obtained joint electricity-carbon prices $\lambda_i^{\text{pro}}$ and $\lambda_j^{\text{pro}}$ are provided in Fig. \ref{fig:price}. With these prices, we can easily check that the optimal dispatch strategies obtained by \eqref{eq:social-optimum} coincide with the optimal solutions of \eqref{eq:strategic-gen} and \eqref{eq:strategic-load} given the obtained joint electricity-carbon prices. This verifies Proposition \ref{prop-3}.

\subsubsection{Performance Comparison}
\begin{figure}[ht]
  \centering
\includegraphics[width=0.5\textwidth]{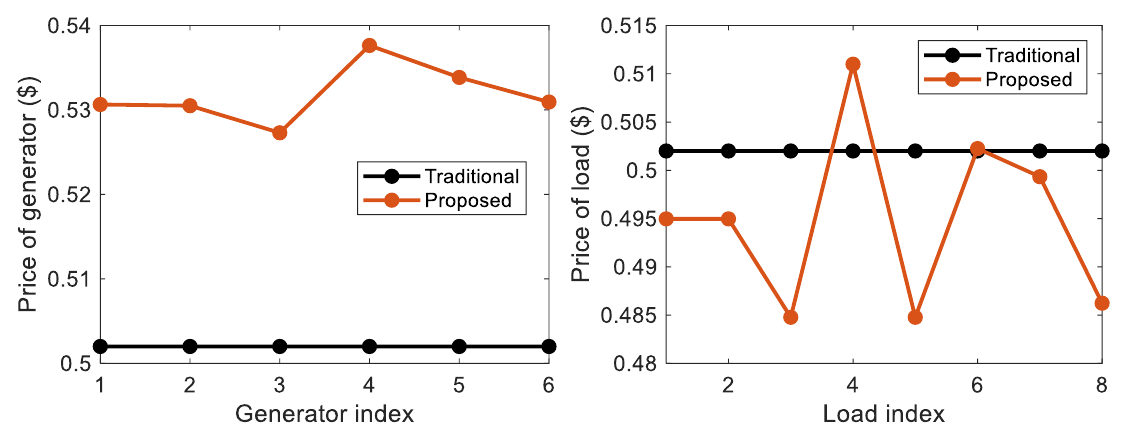}\\
  \caption{Marekt prices of generators and loads under the traditional and proposed markets.}\label{fig:price}
\end{figure}

We first compare the proposed market and the traditional electricity market without carbon considerations. The market prices of generators and loads under these two markets are shown in Fig. \ref{fig:price}. As we can see, the traditional electricity market has a unified price for generators and loads when there is no network congestion. By contrast, in the proposed market, the prices at different buses vary. The generator with a lower carbon intensity has a higher price $\lambda_i^{\text{pro}}$ so that it can earn more from the market and has the incentive to produce more. 

\begin{figure}[ht]
  \centering
\includegraphics[width=0.45\textwidth]{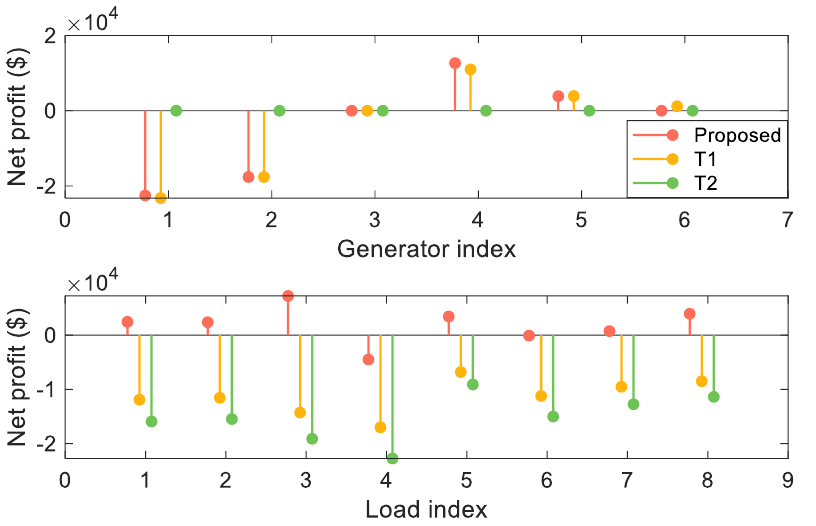}\\
  \caption{Profit change under the Proposed, T1, and T2 methods.}\label{fig:netprofit}
\end{figure}

We further compare the traditional electricity market (without considering carbon emission), T1, T2, and the proposed method in Table \ref{tab:Proposed-T1-T2}. The net profits of individual generators and loads under the Proposed, T1, and T2 methods minus those under the traditional electricity market are recorded in Fig. \ref{fig:netprofit}. 
From the figure, we can find that both T1 and T2 lead to a significant reduction in the net profits of loads, making them reluctant to support the low-carbon transition. By contrast, the net profits of loads change a little under the proposed method. This can also be observed from Table \ref{tab:Proposed-T1-T2}: compared to the traditional electricity market, the total net profit of loads, under T1 and T2, decreases by 12.60\% and 16.86\%, respectively, while that of the proposed method increases by 2.17\%. 
Both the traditional and the proposed markets can achieve budget balance, with subsidies equal to zero. Both T1 and T2 require negative subsidies (\$-107,520 and \$-121,520), meaning that they are generating profits. However, these profits are at the cost of much higher load payments. Further, we also care about whether the prices obtained by different methods can motivate generators and loads to act toward the carbon-aware social optimum. Looking at the last column of Table \ref{tab:Proposed-T1-T2}, we can find that the social welfare under T1 and the proposed method equal the optimal objective value of \eqref{eq:social-optimum}, implying that they have the property of dispatch-following incentive compatibility. Compared to them, T2 results in a lower social welfare, which is not carbon-aware social optimal.

\begin{table*}[t]
\renewcommand{\arraystretch}{1.3}
\renewcommand{\tabcolsep}{0.5em}
\centering
\caption{Comparison of Different Pricing Methods (Unit: \$)}
\label{tab:Proposed-T1-T2}
\begin{tabular}{ccccccccc} 
\hline
Methods & Generator net profit & Load net profit & Generator revenue & Load payment &  Carbon tax & Subsidy & Social Welfare \\
\hline
Traditional & 60,550	& 720,760
 & 1,340,340	& 1,340,340	& 0 & 0   & 659,790\\
T1 &  35,850	& 629,980 & 1,431,120 & 1,431,120 & 107,520 & -107,520 & 665,830
 \\
T2 & 60,550 & 599,240 
& 1,340,340 & 1,461,860 & 0 & -121,520 & 659,790\\
Proposed & 36,946 & 736,400  & 1,421,658 & 1,324,696 & 96,962 & 0 & 665,830\\
\hline
\end{tabular}
\end{table*}

\subsection{IEEE 39-bus system}
\begin{figure}[ht]
  \centering
\includegraphics[width=0.5\textwidth]{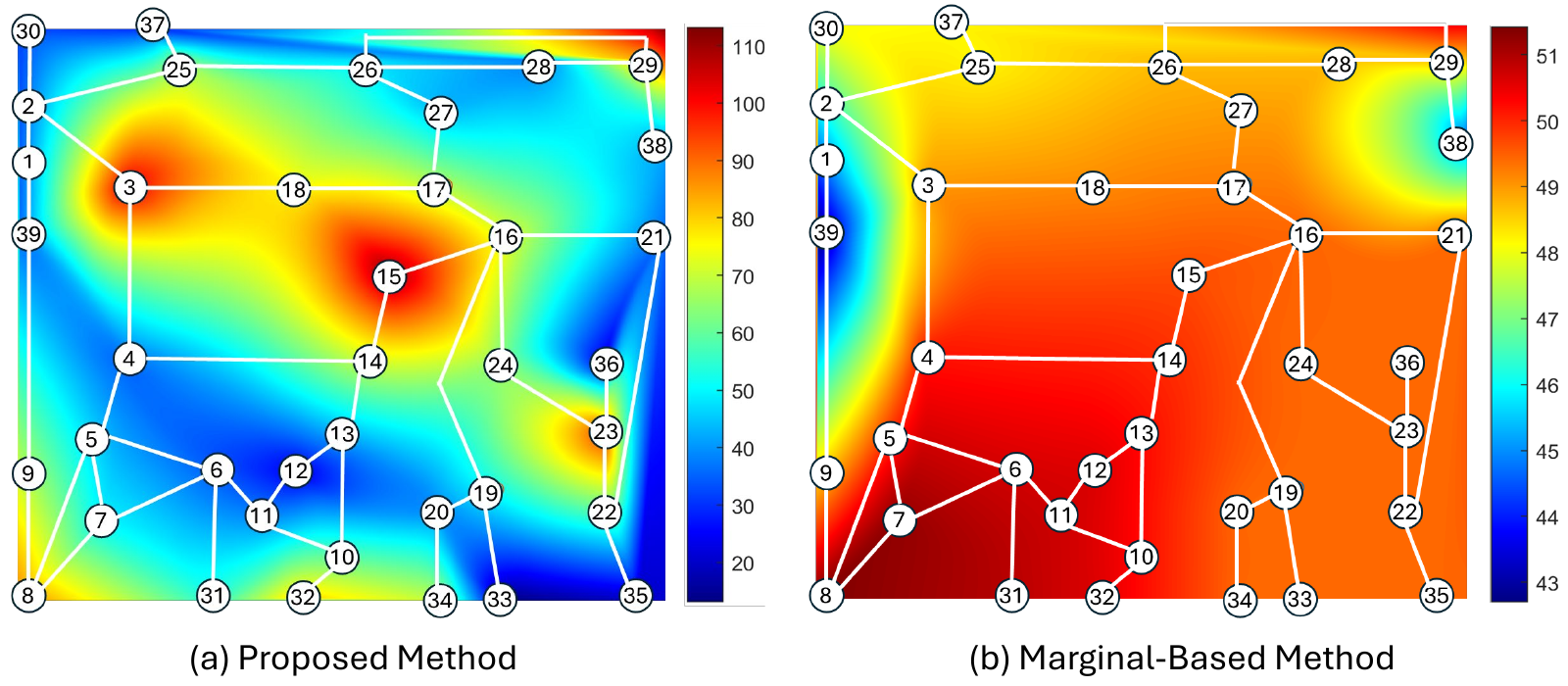}\\
  \caption{Joint electricity-carbon prices at different buses under the proposed and the marginal-based methods.}\label{fig:price-compare}
\end{figure}

\begin{figure}[ht]
  \centering
\includegraphics[width=0.4\textwidth]{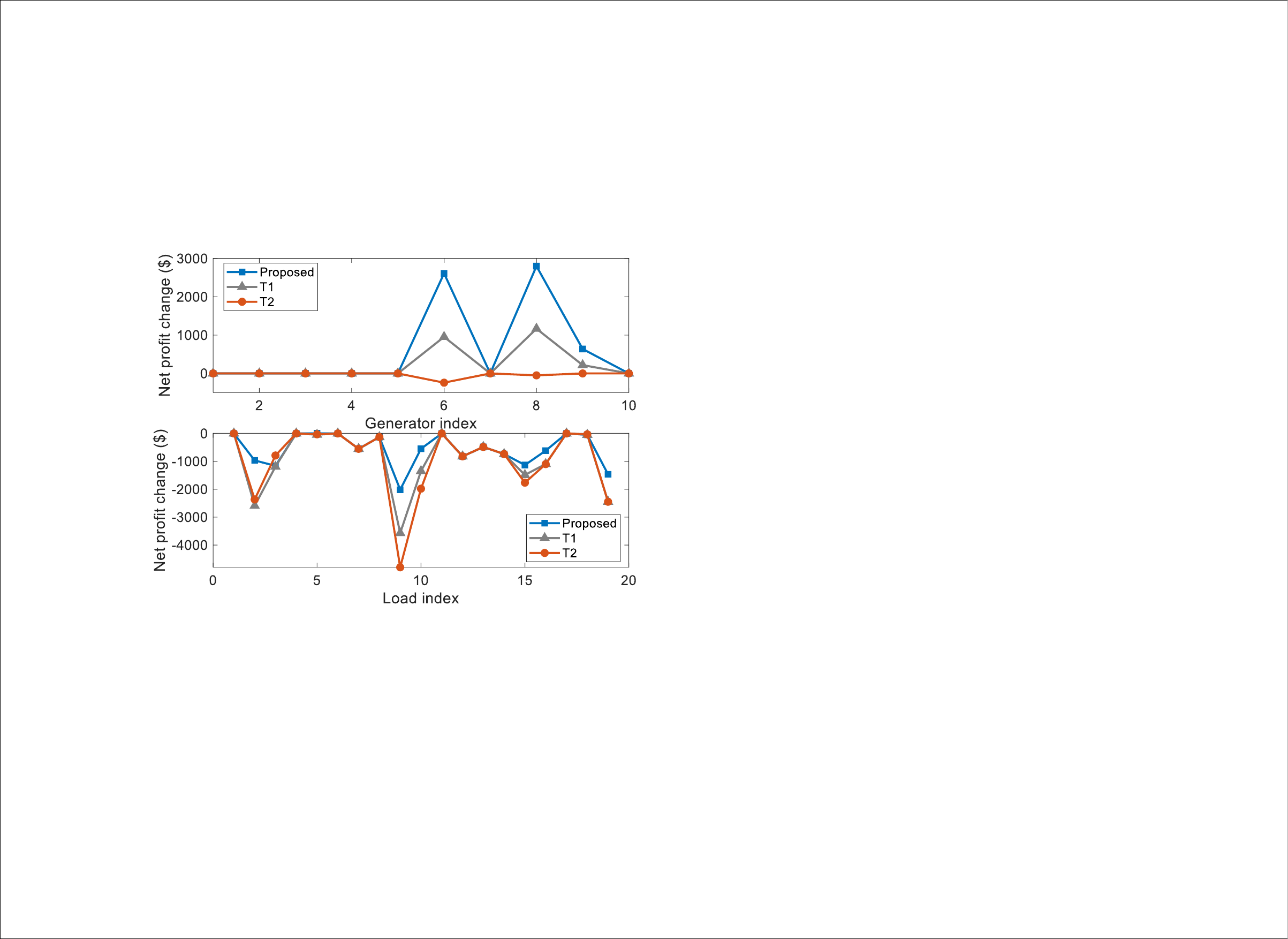}\\
  \caption{Net profit changes of generators and loads under three different methods.}\label{fig:profit-compare}
\end{figure}

We use the IEEE 39-bus system to analyze the case with power network congestion. First, we compare the joint electricity-carbon prices obtained by the proposed method and the marginal-based method (T1 Method), as shown in Fig. \ref{fig:price-compare}. We can find that the prices by the proposed method can better incentivize low-carbon power system operation. Specifically, the generator located at bus 32 has the lowest carbon intensity (0.032 kgCO$_2$/kWh). Correspondingly, its selling price in Fig. \ref{fig:price-compare}(a) is relatively high to encourage it to produce more. The generator at bus 36 has the highest carbon intensity (0.76 kgCO$_2$/kWh) and its selling price is low, which reduces its production. By contrast, the prices at different locations under the T1 Method vary little. The net profit changes of generators and loads under the proposed, T1, and T2 methods are provided in Fig. \ref{fig:profit-compare}. While all methods reduce consumer net profits, the proposed method results in the smallest reduction. Therefore, the consumers are less reluctant to the transition to joint electricity-carbon prices.

\begin{table}[t]
\renewcommand{\arraystretch}{1.3}
\renewcommand{\tabcolsep}{0.5em}
\centering
\caption{Total welfare \eqref{eq:social-optimum-1} under different methods}
\label{tab:flowchange}
\begin{tabular}{cccccc} 
\hline
 & Traditional & T1 & T2 & Proposed & Social Optimum \\
\hline
$0.8F_l$ & -732.08 & 8252.51 & 3062.75 & 8252.51 & 8252.51 \\ 
$0.9F_l$ & -459.54 & 8599.04 & 3202.87 & 8599.04 &  8599.04 \\
$F_l$ & 387.72 & 8819.07 &  3075.96 & 8819.07 & 8819.07 \\
$1.1F_l$ & -134.69 & 8967.12 & 2221.33 & 8967.12 & 8967.12\\
$1.2F_l$ & -1487.24 & 9041.53 & 2484.55 & 9041.53 & 9041.53\\
$1.3F_l$ & -2694.38 & 9115.94 & 2032.88 & 9115.94 & 9115.94 \\
\hline
\end{tabular}
\end{table}

Further, we change the power flow limits $F_l,\forall l\in\mathcal{L}$ from 0.8 to 1.3 times of their original values and calculate the social welfare \eqref{eq:social-optimum-1} under the Proposed, T1, and T2 methods. The social welfare under the traditional electricity market without considering carbon emission is also calculated for comparison. The results are recorded in Table \ref{tab:flowchange}. 
The proposed and T1 methods can always achieve carbon-aware social optimum. Moreover, their social welfare increases as the power flow limits are relaxed. This is straightforward because a less stringent power flow limit provides more flexibility to reduce the operation costs and carbon emissions. The social welfare under T2 method is much lower while the traditional electricity market could result in a negative social welfare.  


\subsection{Scalability}
To show the scalability of the proposed method, we apply it to the simple system in Section \ref{sec:simple system}, the IEEE 39-bus system, and the IEEE 118-bus system. The computation time under different power flow limits are shown in Table \ref{tab:time}. The computation times are all less than 1s, showing the scalability of the proposed method.

\begin{table}[h]
\renewcommand{\arraystretch}{1.3}
\renewcommand{\tabcolsep}{0.5em}
\centering
\caption{Computation time (Unit: second)}
\label{tab:time}
\begin{tabular}{ccccccc} 
\hline
 & $0.8F_l$ & $0.9F_l$ & $1.0F_l$  & $1.1F_l$  & $1.2F_l$  &    $1.3F_l$ \\ 
\hline
Simple system &  0.29 & 0.30 &  0.29  &  0.39 & 0.33 & 0.29  \\
IEEE 39-bus & 0.47 & 0.34 & 0.36 & 0.78  &  0.52  &   0.32\\
IEEE 118-bus & 0.65 & 0.71 & 0.50 & 0.57 &  0.56 &  0.54\\
\hline
\end{tabular}
\end{table}


\section{Conclusion}
\label{sec-V}
This paper proposes a novel joint electricity-carbon pricing mechanism based on a primal-dual optimality condition enabled transformation. Several desired properties of the proposed pricing mechanism are proven theoretically, including budget balance, individual rationality, dispatch-following incentive compatibility, and truthful-bidding incentive compatibility. Case studies reveal the following findings:
\begin{itemize}
    \item The proposed pricing mechanism has less impact on the net profits of generators and loads compared to the two existing methods.
    \item Both the proposed and the marginal-based methods can achieve carbon-aware social optimum, but the CEF-based method may not.
    \item  The prices under the proposed mechanism can better incentivize carbon-intensive generators to produce less.
\end{itemize}

In future research, we may extend the proposed pricing mechanism to incorporate bounded rationality of market participants and uncertainties in generation, load, and the system. We may also explore the applications of the proposed pricing mechanism in various scenarios.



\appendices
\makeatletter
\@addtoreset{equation}{section}
\@addtoreset{theorem}{section}
\makeatother
\setcounter{equation}{0}  
\renewcommand{\theequation}{A.\arabic{equation}}
\renewcommand{\thetheorem}{A.\arabic{theorem}}
\section{Proof of Proposition \ref{prop-1}}
\label{appendix-1}
Let $\upsilon_{\delta}^*$ be the collection of all the primal and dual variables of \eqref{eq:market-proposed} at the optimum under $\delta$. For notation conciseness, in the following, we use $\upsilon_{x}^*$, $\eta_{x}^*$, $\tau_{x}^*$, etc. to represent the optimal values of the corresponding variables under $\delta=x$.

First, we consider two extreme cases: 
\begin{itemize}
    \item When $\delta=0$, we have $\Delta \mathcal{S}=\mathcal{S}_3 \ge 0$.
    \item When $\delta=1$, we can easily get an optimal $\upsilon_1^*$ with $\eta_1^*=0$ by letting $\tau_1^*=\lambda_1^*$. As a result, $\Delta \mathcal{S}=\mathcal{S}_2<0$.
\end{itemize}

Before we proceed, we prove two lemmas as follows.

\begin{lemma} \label{lemma-1}
    If $\eta_{\pi}^*=0$ under $\pi \in [0,1)$, then for all $\delta \in [\pi,1]$, we have $\eta_{\delta}^*=0$.
\end{lemma}

\emph{Proof of Lemma \ref{lemma-1}}: Suppose $\upsilon_{\pi}^*$ and $\upsilon_1^*$ are the optimal under $\delta=\pi$ and $\delta=1$, respectively. 
Then, they both meet the primal-dual optimality condition of \eqref{eq:market-proposed}. In particular, this condition constitutes a set of linear constraints that form a convex set. Let 
\begin{align}
\upsilon_{\delta}^*=\frac{1-\delta}{1-\pi} \upsilon_{\pi}^* + \frac{\delta-\pi}{1-\pi} \upsilon_1^*.
\end{align}
Since $\frac{1-\delta}{1-\pi} \in [0,1]$, $\frac{\delta-\pi}{1-\pi} \in [0,1]$, and $\frac{1-\delta}{1-\pi}+\frac{\delta-\pi}{1-\pi}=1$, according to the definition of a convex set, $\upsilon_{\delta}^*$ also satisfies the primal-dual optimality condition of \eqref{eq:market-proposed} under $\delta$. \hfill $\blacksquare$

With Lemma \ref{lemma-1}, we first keep decreasing $\delta$ from 1 until $\eta_{\delta}^*$ becomes non-zero. Denote the last (lowest) $\delta$ with $\eta_{\delta}^*=0$ as $\tilde \delta$. Then, we have
\begin{lemma}\label{lemma-2}
    For $\delta \in (0,\tilde \delta)$, we have $\eta^*_{\delta}=\eta_{0}^* (1-\delta / \tilde\delta )$.
\end{lemma}

\emph{Proof of Lemma \ref{lemma-2}}: 
Similar to the proof of Lemma \ref{lemma-1}, we have that 
\begin{align}
\upsilon^*_{\delta}=\upsilon_{0}^*(1- \delta / \tilde\delta)+ \upsilon_{\tilde \delta}^* \delta / \tilde\delta  \label{eq:proof.lemma.2:1}
\end{align}
satisfies the primal-dual optimality condition of \eqref{eq:market-proposed} under $\delta \in (0,\tilde \delta)$. Second, we need to prove that $\eta_{\delta}^*$ is the smallest among all optimal $\eta^*$ under $\delta \in (0,\tilde \delta)$. This can be proven by contradiction: if there exists $\check \eta_{\delta}^* < \eta_{\delta}^*$, let
\begin{align}
    \hat \delta := \delta / (1- {\check \eta_{\delta}^*}/{\eta_0^*}) < \delta / (1- {\eta_{\delta}^*}/{\eta_0^*}) =  \tilde \delta.\label{eq:proof.lemma.2:2}
\end{align}
Then using the same property that led to \eqref{eq:proof.lemma.2:1}, we have
\begin{eqnarray}\nonumber
\check \eta_{\delta}^* = \eta_0^* \left(1-{\delta}/{\hat \delta}\right) + \eta^*_{\hat\delta}{\delta}/{\hat \delta}
\end{eqnarray}
which together with \eqref{eq:proof.lemma.2:2} implies $\eta_{\hat\delta}^*=0$. According to Lemma \ref{lemma-1}, we have for all $\delta \in [\hat \delta,1]$, $\eta_{\delta}^*=0$, which contradicts the definition of $\tilde \delta$. This completes the proof of Lemma \ref{lemma-2}. \hfill $\blacksquare$

Based on Lemmas \ref{lemma-1} and \ref{lemma-2}, $\mathcal{S}_3$ decreases with $\delta$ since $\mathcal{S}_3$ equals $\eta_{\delta}^*$ times a nonnegative constant $O^*$, and $\eta_{\delta}^*$ decreases with $\delta$. Therefore, $\Delta \mathcal{S}=\mathcal{S}_2+\mathcal{S}_3$ decreases with $\delta$. Moreover, $\Delta \mathcal{S}$ is continuous in $\delta$ with $\Delta \mathcal{S}(\delta=0) \ge 0$ and $\Delta \mathcal{S}(\delta=1) <0$, so there exists $\delta \in [0,1)$ that makes $\Delta \mathcal{S}=0$. 

In addition, if there is no network congestion, i.e., none constraint in  \eqref{eq:social-optimum-5} is binding, we have $\mathcal{S}_1=0$ since $\underline{\alpha}_l^*=\overline{\alpha}_l^*=0$, $\forall l\in\mathcal{L}$. Hence $\mathcal{S}=0$, i.e., the overall budget balance is achieved. \hfill $\blacksquare$

\setcounter{equation}{0}  
\renewcommand{\theequation}{B.\arabic{equation}}
\renewcommand{\thetheorem}{B.\arabic{theorem}}
\section{Proof of Proposition \ref{prop-2}}
\label{appendix-2}
Under the proposed pricing mechanism, each generator $i \in \mathcal{I}$ solves \eqref{eq:strategic-gen}. Obviously, $p_i=0$ is a feasible solution to \eqref{eq:strategic-gen} with its objective value zero. Since problem \eqref{eq:strategic-gen} is a maximization problem, we have its optimal objective value $\text{NP}_i^* \ge 0$. Similarly, each load $j\in\mathcal{J}$ solves \eqref{eq:strategic-load} with $d_j=0$ a feasible solution. Since problem \eqref{eq:strategic-load} is a maximization problem, we have $\text{NU}_j^*\ge 0$. \hfill $\blacksquare$

\setcounter{equation}{0}  
\renewcommand{\theequation}{C.\arabic{equation}}
\renewcommand{\thetheorem}{C.\arabic{theorem}}
\section{Proof of Proposition \ref{prop-3}}
\label{appendix-3}
Given the realized prices $\lambda_i^{\text{pro}},\forall i\in\mathcal{I}$ and $\lambda_j^{\text{pro}},\forall j\in\mathcal{J}$, the net profit- or utility-maximizing generators and loads solve \eqref{eq:strategic-gen} and \eqref{eq:strategic-load} to determine their dispatch $\tilde{p}_i,\forall i\in\mathcal{I}$ and $\tilde{d}_j,\forall j\in\mathcal{J}$, respectively. Due to separability of the Lagrangian $\mathscr{L}$, $(\tilde p, \tilde q)$ is the optimal solution of \eqref{eq:Lagrangian-2}. Since Slater's condition holds, problems \eqref{eq:Lagrangian-2} and \eqref{eq:market-proposed} have the same optimal solution. 
Furthermore, constraints \eqref{eq:market-proposed-2}--\eqref{eq:market-proposed-9} are the primal-dual optimality condition of problem \eqref{eq:social-optimum}. Thus, problem \eqref{eq:market-proposed} has a unique feasible (also optimal) solution $(p^*,d^*)$ which is the optimal solution of \eqref{eq:social-optimum}. Therefore, $(\tilde p, \tilde q)$ is also the optimal solution of \eqref{eq:social-optimum}, which equals $(p^*,q^*)$. \hfill $\blacksquare$

\setcounter{equation}{0}  
\renewcommand{\theequation}{D.\arabic{equation}}
\renewcommand{\thetheorem}{D.\arabic{theorem}}
\section{Proof of Proposition \ref{prop-4}}
\label{appendix-4}

Let $\theta$ denote the bidding parameters of a pricing-taking agent. Let $\theta^*$ denote the true value of the parameters. Taking a generator $i\in\mathcal{I}$ for example, the optimization model to determine its optimal bidding strategy is
\bsq \label{eq:price-taking}
\begin{align}
    \max_{\theta} ~ & - (c_i + \delta \kappa e_i) p_i^*(\theta) + \lambda_i^{\text{pro}} p_i^*(\theta) \label{eq:price-taking-1} \\~
     \text{s.t.}~ & 0 \le p_i^*(\theta) \le \overline{P}_i \label{eq:price-taking-2}
\end{align}
\esq
where $\lambda_i^{\text{pro}}$ in \eqref{eq:price-taking} is taken as a given constant  under the price-taking assumption. $p_i^*(\theta)$ is the optimal solution of \eqref{eq:social-optimum} with the bidding parameters $\theta$. Moreover, $p_i^*(\theta)$ satisfies \eqref{eq:price-taking-2}.

According to Proposition \ref{prop-3}, we have
\begin{align}
    ~ & - (c_i +\delta \kappa e_i) p_i^*(\theta^*) + \lambda_i^{\text{pro}} p_i^*(\theta^*) \nonumber\\
    \ge ~ & - (c_i +\delta \kappa e_i) p_i + \lambda_i^{\text{pro}} p_i, ~\forall p_i \in [0,\overline{P}_i]
\end{align}
Therefore, bidding truthfully (bidding $\theta^*$) is the optimal for a price-taking generator. The case for loads can be proven in a similar way. \hfill $\blacksquare$

\bibliographystyle{IEEEtran}
\bibliography{PaperRef}

\end{document}